\documentclass{amsart}

\usepackage{amsmath,amsthm,amssymb}
\usepackage{amsfonts}
\usepackage{rotating}
\usepackage{euscript}
\usepackage{pst-node}
\usepackage{epsfig}

\def\be{\begin{equation}}
\def\ee{\end{equation}}

\def\C{{\mathbb C}} 
\def\N{{\mathbb N}} 
\def\P{{\mathbb P}}

\def\R{{\mathbb R}}

\def\f{\EuScript}

\def\phi{{\varphi}}
\def\v{{\varepsilon}} 
\def\deg{{\rm deg\,}}

\def\cos{{\rm cos\,}} 
 
\def\GCD{{\rm GCD }}

\def\tilde{\widetilde}

\def\bp{\begin{proposition}}
\def\ep{\end{proposition}}

\def\bt{\begin{theorem}}
\def\et{\end{theorem}}
\def\br{\begin{remark}}
\def\er{\end{remark}}
\def\be{\begin{equation}}
\def\bee{\begin{equation*}}
\def\la{\label}

\def\ee{\end{equation}}
\def\eee{\end{equation*}}
\def\bl{\begin{lemma}}
\def\el{\end{lemma}}
\def\bc{\begin{corollary}}
\def\ec{\end{corollary}}
\def\pr{\noindent{\it Proof. }}

\def\bd{\begin{definition}}
\def\ed{\end{definition}}
\newcommand{\cal}{\mathcal}
\input epsf.sty

\newtheorem{theorem}{Theorem}[section]
\newtheorem{lemma}{Lemma}[section]
\newtheorem{definition}{Definition}[section]
\newtheorem{corollary}{Corollary}[section]
\newtheorem{proposition}{Proposition}[section]
\newtheorem{remark}{Remark}[section]

\begin{document}
\title{On decompositions of trigonometric polynomials}
\author{F. Pakovich}
\date{}
\begin{abstract} 
Let $\R_t[\theta]$ be the ring  generated over $\R$ by 
$\cos \theta$ and $\sin \theta$, and $\R_t(\theta)$ be its quotient  field.
In this paper we study the ways in which an element $p$ of $\R_t[\theta]$
can be decomposed into a composition of functions of the form 
$p=R\circ q$, where $R\in \R(x)$ and $q\in \R_t(\theta)$.
In particular, we describe all possible solutions of the functional equation   
$R_1\circ q_1=R_2\circ q_2,$ where $R_1, R_2 \in \R[x]$ and $q_1,q_2\in\R_t[\theta].$

\end{abstract}
\maketitle

\section{Introduction} 
Let $P$ be a polynomial with complex coefficients. Any representation of $P$ in the form $P=P_1\circ W_1$, where $P_1$ and $W_1$ are 
polynomials of degree greater than one and the symbol $\circ$ denotes the superposition of functions, is called a decomposition of $P.$  
The problem of description of all possible decompositions of a polynomial  
naturally leads to the functional equation 
\be \label{ma} P_1\circ W_1=P_2\circ W_2,\ee
where $P_1,W_1,P_2,W_2$ are polynomials, for the first time studied by Ritt in the paper \cite{r2}. 
In particular, the results of \cite{r2} imply that in a certain sense  all polynomial solutions of \eqref{ma} reduce  
either to the solutions 
$$z^n \circ z^rR(z^n)=z^rR^n(z) \circ z^n,$$
where $R$ is a polynomial, and $r\geq 0,$ $n\geq 1,$ or to the solutions 
\be \label{gfd} T_n \circ T_m= T_m \circ T_n,\ee
where $T_n, T_m$ are the Chebyshev polynomial.

Functional equation \eqref{ma} is closely
related to the so-called ``polynomial moment problem"
which asks to describe complex polynomials $P,Q$ such
that the equalities 
\be \label{4} \int_{0}^1 P^idQ=0, \ \ \ i\geq 0,\ee hold.
Indeed, it is easy to see using the change $z\rightarrow W(z)$ that \eqref{4} is satisfied whenever there exist 
polynomials $\widetilde P,$ $\widetilde Q,$ and $W$ such that
\be \label{2}
P=\widetilde P\circ W, \ \ \ 
\ Q=\widetilde Q\circ W, \ \ \ W(0)=W(1). 
\ee
Furthermore, it was shown in \cite{pm} that if polynomials $P,$  $Q$ satisfy \eqref{4}, then there exist polynomials $Q_j$ such that $Q=\sum_j Q_j$ and the equalities 
\be \label{cc}
P=\tilde P_j\circ W_j, \ \ \
Q_j=\tilde Q_j\circ W_j, \ \ \ W_j(0)=W_j(1)
\ee hold 
for some polynomials $\tilde P_j, \tilde Q_j, W_j$. Thus, the most interesting solutions of the polynomial moment problem 
arise from polynomials having ``multiple'' decompositions 
\be \label{sol} P=\widetilde P_1\circ W_1=\widetilde  P_2\circ W_2=\dots=\widetilde P_s\circ W_s.\ee
Polynomial solutions of \eqref{sol} were described in the paper \cite{pakk}, where the correspon\-ding generalization of the result of Ritt about solutions of \eqref{ma} 
was obtained.

The polynomial moment problem naturally appears in the study of the center problem for the Abel differential equation with polynomial coefficients which is a simplified analog 
of the center problem  for the Abel differential equation whose coefficients are 
trigonometric polynomials over $\R$  (see e. g. the recent papers \cite{bry}, \cite{bpy} and the bibliography therein). In its turn, the last problem
is closely related to the classical center-focus problem of Poincar\'e 
(\cite{cher}).
In the same way as the center problem for the Abel equation with polynomial coefficients leads to the polynomial moment problem,  
the center problem for the Abel equation with trigonometric coefficients leads to the following ``trigonometric moment problem''.
Let $$p=p(\cos \theta,\sin \theta),\ \ \ q=q(\cos\theta ,\sin\theta )$$ be trigonometric polynomials over $\R$, that is elements of the ring $\R_t[\theta]$ generated over $\R$ by the functions 
$\cos \theta$, $\sin \theta$.
What are conditions implying that the equalities 
\be \label{1} \int_0^{2\pi }p^idq=0, \ \ \ i\geq 0,\ee
hold ? Like to the case of the polynomial moment problem one can consider a complexified version of this problem 
(see  
\cite{ppre}, \cite{ppz}, \cite{abc}). However, examples constructed in \cite{ppz}, \cite{abc} suggest that in the trigonometric case the complex version of the problem may be much more complicated than the real one.

Again, a natural sufficient condition for \eqref{1} to be satisfied is related with compositional properties of $p$ and $q$. 
Namely, it is easy to see that if there exist $P, Q\in \R[x]$ and $w\in \R_t[\theta]$ such that  
\be \label{c} p=P\circ w, \ \ \ \ q=Q\circ w, \ee
then \eqref{1} hold. Furthermore, if for given $p$ there exist several such $q$ (with different $w$),
then \eqref{1} obviously holds for their sum.
Thus, the trigonometric moment problem leads to the problem of description of solutions 
of the equation 
\be \label{maii} P_1\circ w_1=P_2\circ w_2,\ee where $w_1,w_2\in \R_t[\theta]$ and $P_1, P_2 \in \R[x],$ and
the main goal of this paper is to provide such a description. Notice that, besides of its relation with the trigonometric moment problem,
functional equation \eqref{maii} 
seems to be 
interesting by itself.  
In particular, it contains among its solutions  the most  
known  trigonometric identity \be \label{osn}\sin^2 \theta=1-\cos^2 \theta.\ee
Besides, the problem of description of solutions of \eqref{maii} absorbs the problem of description of polynomial solutions of 
\eqref{ma} over $\R$ since for any 
polynomial solution of \eqref{ma} and any $w\in \R_t[\theta]$
we obtain a solution of  \eqref{maii} setting 
$$w_1=W_1\circ w, \ \ \  w_2=W_2\circ w.$$

Observe that if $P_1,P_2,w_1,w_2$ is a solution of \eqref{maii}, then  for any $k\in \N$ and $b\in \R$ 
we obtain another solution $P_1,P_2,\widetilde w_1,\widetilde w_2$   setting
$$ \widetilde w_1(\theta) = w_1(k\theta+b) , \ \ \ \widetilde w_2(\theta) = w_2(k\theta+b). $$ 
Further, if $P_1,P_2,w_1,w_2$ is a solution of \eqref{maii}, then 
for any $U\in \R[t]$ we obtain another solution $\widetilde P_1,\widetilde P_2,w_1,w_2$ 
setting
$$\widetilde P_1=U\circ P_1, \ \ \  \widetilde P_2=U\circ P_2.$$

Let $p$ be an element of $\R_t[\theta]$ or $\R[x],$ and
$p=P_1\circ w_1$ and $p=\widetilde P_1\circ \widetilde w_1$ be two decompositions of 
$p$,  such that $P_1,\widetilde P_1\in \R[x]$ and $w_1,\widetilde w_1 \in \R_t[\theta]$ or $w_1,\widetilde w_1 \in \R[x]$.  
We will call  
these decompositions  equivalent, and use the notation $P_1\circ w_1\sim \widetilde P_1\circ \widetilde w_1$,
if there exists $\mu\in \R[x]$ of degree one such that 
$$\widetilde P_1=P_1\circ \mu, \ \ \ \widetilde w_1=\mu^{-1}\circ w_1.$$

Under the above notation our  main result about solutions of \eqref{maii} may be formulated as follows.

\bt Assume that 
$P_1, P_2 \in \R[x]\setminus \R$ and  $w_1,w_2\in \R_t[\theta]\setminus \R$ satisfy the equality 
$$P_1\circ w_1=P_2\circ w_2.$$ 
Then,  up to a possible replacement of $P_1$ by $P_2$ and $w_1$ by $w_2$, one  of the following conditions holds.
\vskip 0.2cm
\noindent  1. \ \  
There exist $U,\widetilde P_1,\widetilde P_2,W_1,W_2\in \R[x]$ and $\widetilde w \in \R_t[\theta]$ such that 
$$ P_1=U\circ \widetilde P_1, \ \ \ P_2=U\circ \widetilde P_2, \ \ \ 
w_1=W_1\circ \widetilde w, \ \ \  w_2=W_2\circ\widetilde w ,\ \ \ \widetilde P_1\circ W_1=\widetilde P_2\circ W_2,$$ 
and either
$$
\widetilde P_1\circ W_1\sim  z^n \circ z^rR(z^n),
  \ \ \ \ \ \ \widetilde P_2\circ W_2\sim  z^rR^n(z) \circ z^n,\leqno\ \ \ a)$$	
where $R\in \R[x]$, $r\geq 0,$ $n\geq 1,$ and $\GCD(n,r)=1$, 
or
$$
\widetilde P_1\circ W_1\sim T_{n} \circ T_m, \ \ \ \ \ \ 
\widetilde P_2\circ W_2\sim T_m\circ T_n, \leqno \ \ \ b)$$
where $T_{n}$ and $T_{m}$ are the Chebyshev polynomials, $m,n\geq 1,$ and  $\GCD(n,m)=1$.
\vskip 0.2cm
\noindent 2.\ \ There exist 
$U, \widetilde P_1, \widetilde P_2\in \R[x],$ $\widetilde w_1,$ $\widetilde w_2\in \R_t[\theta],$ and a polynomial $W(\theta)=k\theta +b,$ where $k\in \N,$ $b\in \R$, such that 
$$ P_1=U\circ \widetilde P_1, \ \ \ P_2=U\circ \widetilde P_2, \ \ \ 
w_1=\widetilde w_1\circ W, \ \ \  w_2=\widetilde w_2\circ W,\ \ \ \widetilde P_1\circ \widetilde w_1=\widetilde P_2\circ \widetilde w_2,$$ 
and either
$$\widetilde P_1\circ \widetilde w_1\sim  z^2 \circ \,
\cos \theta\, S(\sin \theta),  \ \ \ \ \ \ \widetilde P_2\circ \widetilde w_2\sim  (1-z^2)\, S^2(z)\circ \sin\theta ,\leqno \ \ \ a)$$
where $S\in \R[x]$, or
$$\widetilde P_1\circ \widetilde w_1\sim -T_{nl} \circ \cos\left( \frac{(2s+1)\pi}{nl}+m\theta\right), \ \ \ \ \ \ 
\widetilde P_2\circ \widetilde w_2\sim T_{ml} \circ \cos(n\theta),\leqno \ \ \ b)$$    
where $T_{nl},T_{ml}$ are the Chebyshev polynomials, $m,n\geq 1,$  $l>1$,  $0\leq s < nl,$ and $\GCD(n,m)=1$.
\et

Notice that solutions of types 1, a) and 1, b) reduce to polynomial solutions of \eqref{ma}, while
solutions of type 2, a) generalize identity \eqref{osn}. Further, solutions 
of type 2, b) can be considered as a generalization of the identity $$T_n\circ \cos m\theta= T_m\circ \cos n\theta,$$ although 
this  identity itself is  an example of  a solution of type 1, b) since $$\cos m\theta=T_m\circ \cos\theta, \ \ \ 
\cos n\theta=T_n\circ \cos\theta.$$

Our approach to functional equation \eqref{maii} relies on the isomorphism
$$\phi:\, \ \cos \theta\rightarrow \left(\frac{z+1/z}{2}\right), \ 
\sin \theta\rightarrow \left(\frac{z-1/z}{2i}\right),$$ 
between the ring $\R_t[\theta]$ and  a subring  of the ring $\C[z,1/z]$ of complex Laurent polynomials.
Clearly, any decomposition $p=P\circ w$ of $p\in \R_t[\theta]$, where $P\in \R[x]$ and $w\in \R_t[\theta]$, or more generally where  
$P\in \R(x)$ and $w$ is contained in the quotient field $\R_t(\theta)$ of $\R_t[\theta]$, descends to a decomposition   
$\phi(p)=P\circ \phi(w)$ of $\phi(p)$, making it possible to use results of \cite{pak}, about decompositions of Laurent polynomials into compositions of rational functions, 
for the study of decompositions of trigonometric polynomials.

The paper is organized as follows. In the second section we recall some basic facts about decompositions of Laurent polynomials and prove their analogues for decompositions in $\R_t[\theta]$.
We also
show (Corollary \ref{c1}) that for $p\in \R_t[\theta]$ any equivalence class of decompositions of $\phi(p)\in\C[z,1/z]$ into a composition of rational functions over $\C$ contains a representative which lifts to a decomposition $p=P\circ w$, where  
$P\in \R(x)$ and $w\in \R_t(\theta)$. This result shows that the decomposition theory for $\R_t[\theta]$ is ``isomorphic"  
to the decomposition theory for a certain subclass of complex Laurent polynomials, 
and permits to deduce results about decompositions in $\R_t[\theta]$
from the ones in $\C[z,1/z]$.
In the third section, basing on the results of the second section and results of \cite{pak} about  decompositions of Laurent polynomial, we prove
Theorem 1.1.

\section{Decompositions in $\R_t[\theta]$ and in $\C[z,1/z]$}
It is well known that $\R_t[\theta]$ is isomorphic to 
a subring of the field $\R(x)$, where the isomorphism $\psi:\,\R_t[\theta]\rightarrow\R(x)$ is defined by the formulas 
\be \label{t0} \psi(\sin \theta)=\frac{2x}{1+x^2}, \ \ \ \psi(\cos \theta)=\frac{1-x^2}{1+x^2}.\ee 
Furthermore, the isomorphism $\psi$ extends to an isomorphism between $\R_t(\theta)$ and $\R(x)$ which maps the generator $\tan(\theta/2)$  of $\R_t(\theta)$ to the generator $x$ of $\C(x),$ 
$$x=\psi\left(\frac{\sin \theta}{1+\cos\theta}\right)=\psi\left(\tan(\theta/2)\right).$$ 
In particular, this implies  by the
L\"uroth theorem that any subfield $k$ of $\R_t(\theta)$ has the form $k=\R(b)$ for some $b\in \R_t(\theta).$ 
In this paper however
instead of the isomorphism $\psi$ 
we will use the 
isomorphism $\phi$ between  the ring $\R_t[\theta]$ and  a subring  of the ring $\C[z,1/z]$ of complex Laurent polynomials,
defined by the formulas
\be \la{t} \phi(\cos \theta)= \frac{z+1/z}{2}, \ \ \ 
\phi(\sin \theta)=\frac{z-1/z}{2i},\ee 
which seems to be more useful for the study of compositional properties of  $\R_t[\theta]$. 

For brevity,  we will denote the ring $\C[z,1/z]$ by ${\cal L}[z]$ and  the image of 
$\R_t[\theta]$ in  ${\cal L}[z]$ under the isomorphism $\phi$ by ${\cal L_{\R}}[z]$. It is easy to see that ${\cal L_{\R}}[z]$ consists of Laurent polynomials 
$L$ such that $\bar L(1/z)=L(z),$ where $\bar L$ denotes the Laurent polynomial obtained from $L$ by the complex conjugation of all its coefficients. 
Clearly, the isomorphism $\phi$ extends to an isomorphism between  $\R_t(\theta)$  and ${\cal L_{\R}}(z)$, where ${\cal L_{\R}}(z)$
consists of rational functions $R$ satisfying the equality 
$\bar R(1/z)=R(z).$

Any decomposition $p=P\circ w$, where $p\in \R_t[\theta]$,
$P\in \R(x),$ and $w\in \R_t(\theta)$, obviously descends to a decomposition   
$\phi(p)=P\circ \phi(w)$, where $\phi(p)\in {\cal L_{\R}}[z]$ and $\phi(w)\in {\cal L_{\R}}(z)$. 
However, it is clear that $L=\phi(p)$ may have decompositions $L=A\circ B,$ where $A,B\in \C(z),$ such that the coefficients of $A$ are not real and 
$B$ is not contained in ${\cal L_{\R}}(z)$. In this context the following simple lemma is useful.

\bl \label{l0} Let $L\in {\cal L_{\R}}(z)\setminus \R$ and let
$L=A\circ B$ be a decomposition of $L$ into a composition of rational functions $A,B\in \C(z).$ 
Then the inclusion $B\in {\cal L_{\R}}(z)$ implies the inclusion $A\in \R(x)$.
\el
\pr Indeed, since $L,B\in {\cal L_{\R}}(z)$, we have:
$$A\circ B=\bar A\circ \bar B\circ 1/z=\bar A\circ B,$$ implying that  $\bar A=A.$ \qed

\vskip 0.2cm

We will call a Laurent polynomial $L$ proper if $L$ is neither a polynomial in $z$, nor a polynomial in $1/z,$ or in other words 
if $L$ has exactly two poles.
The lemma below is a starting point of the decomposition theory of Laurent polynomials (see \cite{pak}).

\bl \label{l1} Let $L=P\circ W$ be a decomposition of  $L\in {\cal L}[z]\setminus \C$ into a composition of rational functions $P,W\in \C(z).$ 
Then there exists $\mu\in \C(z)$ of degree one
such that either $P\circ \mu $ is a polynomial and $\mu^{-1}\circ W$ is a Laurent polynomial, or
$P\circ \mu $ is a Laurent polynomial and $\mu^{-1}\circ W=z^d$, $d\geq  1.$
\el

\pr Indeed, it follows easily from 
$$L^{-1}\{\infty \}=W^{-1}\{P^{-1}\{\infty \}\}\subseteq \{0,\infty\}$$ that either $P^{-1}\{\infty \}$ consists of a single point $a\in \C\P^1$ and $W^{-1}\{a \}\subseteq\{0,\infty\},$
or $P^{-1}\{\infty \}$ consists of two points $a,b\in \C\P^1$ and $W^{-1}\{a,b \}=\{0,\infty\}.$
In the first case 
there exists a rational function 
$\mu\in \C(z)$ of degree one such that $P\circ \mu $ is a polynomial and $\mu^{-1}\circ W$ is a Laurent polynomial (which is proper if 
and only if $L$ is proper). In the second case there exists $\mu\in \C(z)$ of degree one such that $P\circ \mu $ is a proper Laurent polynomial and $\mu^{-1}\circ W=z^d,$ $d\geq 1$.
\qed
\vskip 0.2cm

The following statement is a ``trigonometric'' analogue of Lemma \ref{l1} and  is equivalent to 
Proposition 21 of \cite{ggl}  and to Theorem 5 of \cite{ga1}. 
Notice however that the proofs given in \cite{ggl}, \cite{ga1}
are much more complicated than the proof given below. The idea to relate decompositions in $\R_t[\theta]$ with decompositions in ${\cal L}[z]$
was proposed in the paper \cite{pakov}, and the proof given below essentially coincides with the proof of Lemma 2.2 in \cite{pakov}. 

\bl \label{l2} Let $p=P\circ w$ be a decomposition of $p\in \R_t[\theta]\setminus \R$ into a composition of $P\in \R(x)$ and $w\in \R_t(\theta).$ Then there exists a rational function 
$\mu\in \R(x)$ of degree one  such that either $P\circ \mu \in \R[x]$ and $\mu^{-1}\circ w\in \R_t[\theta]$, or
$P\circ \mu\in \R(x)$ and $\mu^{-1}\circ w=\tan(d\theta/2)$, $d\geq  1.$ 
\el 
\pr 
Setting $$L=\phi(p), \ \ \ W=\phi(w)$$ and considering the equality $L=P\circ W$,
we conclude as above that either
\be \label{qaz} P^{-1}\{\infty \}=\{a\} \ \ {\rm and} \ \ W^{-1}\{a \}=\{0,\infty\}\ee for some $a\in \C\P^1,$ or  
\be \label{qaz1} P^{-1}\{\infty \}=\{a,b\} \ \  {\rm and} \ \  W^{-1}\{a,b \}=\{0,\infty\}\ee
for some $a, b\in \C\P^1.$

Assume that \eqref{qaz} holds. Since $P\in \R(x)$, it follows from $P^{-1}\{\infty \}=\{a\}$ that either $a\in \R$, or $a=\infty$ and 
$P\in \R[x]$, $W\in {\cal L_{\R}}[z]$.
In the second case, since 
$\phi$ is an isomorphism between $\R_t[\theta]$ and ${\cal L_{\R}}[z],$ we conclude that 
$w\in \R_t[\theta].$ On the other hand, if $a\in \R$, then setting $\mu=a+1/z$ we see that 
$P\circ \mu \in \R[x]$ and $\mu^{-1}\circ W\in {\cal L}[z]$. Furthermore, since $W\in {\cal L_{\R}}(z)$ and $\mu$ has real coefficients, 
the function $\mu^{-1}\circ W$ is contained in ${\cal L_{\R}}[z]$ implying 
that $\mu^{-1}\circ w\in\R_t[\theta].$

If \eqref{qaz1} holds, then we can modify  $\mu\in \C(z)$ from Lemma \ref{l1} so that 
\be \la{iuy} \mu^{-1}\circ W=\frac{1}{i}\frac{z^d-1}{z^d+1}=\frac{1}{i}\left(\frac{z^{d/2}-z^{-d/2}}{z^{d/2}+z^{-d/2}}\right)=\phi(\tan(d\theta/2)), \ \ \ d\geq 1.\ee
Furthermore, since the functions $\phi(\tan(d\theta/2))$ and $W$ are contained in ${\cal L_{\R}}(z)$, it follows from Lemma \ref{l0} that
$\mu^{-1}\in \R(x)$. Therefore, $P\circ \mu\in \R(x)$. Finally, clearly, 
$\mu^{-1}\circ w=\tan(d\theta/2).$ 
\qed
\vskip 0.2cm

Notice that if $p=P\circ w$ is a decomposition of $p\in \R_t[\theta]$ such that $P\in \R(x)$
and $w=\tan(d\theta/2)$, $d\geq  1,$ then $P$ has the form $$P=\frac{A}{(x^2+1)^k},\ \ \ \ k\geq 1,$$   where $A\in \R[x],$  and 
$\deg A\leq 2k,$ since \eqref{iuy} implies that the function  $\mu^{-1}\circ W$ 
sends $0$ and $\infty$ to $i$ and $-i.$ 
Alternatively, we can observe that $\tan(d\theta/2)$ considered as a function of complex variable takes all the values in $\C\P^1$ distinct from $\pm i$.
Therefore, 
the function $P$ may have poles only at points $\pm i$, since otherwise the composition 
$p=P\circ w$ would  not be an entire function.

\vskip 0.2cm

Two different types of decompositions of Laurent polynomials  appearing in Lem\-ma \ref{l1}
correspond to two different types of 
imprimitivity systems in their monodromy groups
(for more details concerning decompositions of rational functions with two poles we refer the reader to \cite{mp}).
Namely, if $L$ is a Laurent polynomial of degree $n$ 
we may assume that its monodromy group $G$ contains the permutation 
$$h=(1\,2\,\dots\, n_1)(n_1+1\, n_1+2\, \dots\, n_1+n_2),$$ where
$1\leq  n_1 \leq n,$ $0\leq n_2 <n,$  
$n_1+n_2=n$. Furthermore, the equalities $n_1=n,$ $n_2=0$ hold if and only if $L$ is not proper.

Let  $\f E$ be an imprimitivity system  of $G$. Denote by $W_{i,d}^1$ (resp. by $W_{i,d}^2$) a union of numbers from the
segment $[1,n_1]$ (resp.  $[n_1+1,n_1+n_2]$) equal to $i$
by modulo $d$. 
Since $h$ permutes blocks   of $\f E$, it is easy to see that 
either there exists a number $d\vert n$  such that
any block of $\f E$ is equal to $W_{i_1,d}^1\cup W_{i_2,d}^2$ for some $i_1,i_2,$  $1\leq i_1,i_2\leq d,$
or  there exist numbers
$d_1\vert n,d_2\vert n$ 
such that
any block of $\f E$ is equal either to $W_{i_1,d_1}^1$ for some $i_1,$ $1\leq i_1\leq d_1,$ or to
$W_{i_2,d_2}^2$ for some $i_2,$  $1\leq i_2\leq d_2.$ Furthermore, since blocks have the same cardinality, in the second case 
\be \la{rav} n_1/d_1=n_2/d_2. \ee
The imprimitivity systems of the first type correspond to decompositions $L=A\circ B,$ where $A$ s a polynomial and $B$ is a Laurent 
polynomial, while imprimitivity systems of the second type correspond to decompositions $L=A\circ B,$ where $A$ is a proper Laurent polynomial and $B=z^d.$

\vskip 0.1cm

The following result coincides with Lemma 6.3 of \cite{pak}. For the reader convenience we provide below  a self-contained proof.

\bl \la{copo1} Let $A,B\in \C[z]\setminus \C$  and $L_1$, $L_2\in {\cal L}[z]\setminus \C$ satisfy 
\be \la{tyr} A\circ L_1=B\circ L_2.\ee Assume additionally that $\deg A = \deg B$.
Then either there exists a polynomial $w\in \C[z]$ of degree one such that 
\be \label{barsu0} B=A\circ w^{-1}, \ \ \ L_2=w\circ L_1,
\ee 
or there exist polynomials $w_1,w_2\in \C[z]$ of degree one such that 
\be \label{barsu} w_1\circ L_1=\left(z^r+\frac{1}{z^r}\right)\circ (az), \ \ \ w_2\circ L_2=\left(z^r+\frac{1}{z^r}\right)\circ (a\nu  z)\ee 
for some $r\in \N$, $a\in \C$, and a root of unity $\nu $.  
 \el

\pr Let $G$ be the monodromy group of a Laurent polynomial $L$ defined by any of the parts of equality \eqref{tyr}. Then  $G$ has two imprimitivity systems  of the first type $\f E_1$ and $\f E_2$, corresponding to the decompositions in \eqref{tyr}. Furthermore, 
since $\deg A = \deg B$, the blocks of $\f E_1$ and $\f E_2$ have the same cardinality $l=\deg L/\deg A.$ 

If these systems coincide, then equalities \eqref{barsu0} hold for 
some rational function 
$w\in \C(z)$ of degree one which obviously is a polynomial.
On the other hand, if they are different, then it is easy to see that the  imprimitivity system $\f E_1\cap \f E_2$
 belongs to the second type, and has blocks consisting of $r$
elements, where $2r=l.$ In particular, $L$ and $L_1,L_2$ are proper, and 
the equalities 
\be \label{rfv} L_1=\tilde L_1\circ W, \ \ \  L_2=\tilde L_2\circ W, \ee
hold for some rational functions 
$\tilde L_1, \tilde L_2, W$, where $\deg \tilde L_1=\deg \tilde L_2=2.$ 
Applying now Lemma \ref{l1} to equalities \eqref{rfv}  we conclude that
$$L_1=\left(\alpha_0+\alpha_1 z+\frac{\alpha_2}{z}\right)\circ z^r, \ \ \ \ L_2=\left(\beta_0+\beta_1z+\frac{\beta_2}{z}\right)\circ z^r,$$
for some $\alpha_0,\beta_0\in \C,$ and $\alpha_1, \alpha_2,\beta_1, \beta_2\in \C\setminus \{0\}$. Furthermore, equality \eqref{tyr} implies that 
$$L_1=\left(\alpha_0+\alpha_1 z+\frac{\alpha_2}{z}\right)\circ z^r, \ \ \ \ L_2=\left(\beta_0+\alpha_1\nu _1 z+\frac{\alpha_2\nu _2}{z}\right)\circ z^r,$$
for some roots of unity $\nu _1,\nu _2$.
The lemma follows now from the equalities 
$$\alpha_0+\alpha_1 z^r+\frac{\alpha_2}{z^r}=\left(\alpha_0+\frac{\alpha_1 z}{a^r} \right)\circ \left(z^r+\frac{1}{z^r}\right)\circ (az),$$ 
$$\beta_0+\alpha_1\nu_1 z^r+\frac{\alpha_2\nu_2}{z^r}=\left(\beta_0+\frac{\alpha_1\nu_1  z}{a^r\nu^r}\right)\circ \left(z^r+\frac{1}{z^r}\right)\circ (a\nu z),$$ where $a$ and $\nu$ are complex numbers satisfying $a^{2r}=\alpha_1/\alpha_2$ and $\nu^{2r}=\nu_1/\nu_2.$ 
\qed
\vskip 0.2cm

\bl \label{uuii}  Let $L=A\circ L_1$ be a decomposition of $L\in  {\cal L_{\R}}[z]\setminus \R$  into a composition
of $A\in \C[z]$ and $L_1=\sum_{-n}^nc_iz^i \in{\cal L}[z].$  Assume additionally that  
$c_{-n}=1/c_{n}.$ Then the leading coefficient of $A$ is real and $\vert c_n\vert =\vert c_{-n}\vert=1.$
\el

\pr
Let $\alpha$ be the leading coefficient of $A$ and $d=\deg A.$
Since $L \in  {\cal L_{\R}}[z]$, we have $\bar \alpha\bar c_{n}^d=\alpha  c_{-n}^d$ implying  that
\be \label{xyu+}\bar \alpha\bar c_n^d = \alpha/c_n^d.\ee  
Multiplying this equality by its conjugated we obtain the equality $(\bar c_n c_n)^{2d}=1.$
Since $\bar c_n c_n=\vert c_n\vert^2$ is a real positive number, we conclude 
that $c_n\bar c_n=1$ or equivalently that $\vert c_n\vert=1$.
Now \eqref{xyu+} implies that $\bar\alpha=\alpha$. \qed

\bt \label{2.1} Let $L=A\circ L_1$ be a decomposition of $L\in  {\cal L_{\R}}[z]\setminus \R$  into a composition
of $A\in \C[z]$ and $L_1 \in{\cal L}[z].$ Then there exists a polynomial $v\in \C[z]$ of degree one
such that $A\circ v^{-1}\in \R[x]$ and $v\circ L_1\in {\cal L_{\R}}[z]$.
\et   

\pr Since $L$ belongs to $\in  {\cal L_{\R}}[z]$, the equality 
$$A\circ L_1=\bar A \circ \bar L_1\circ 1/z$$ holds. Applying to this equality Lemma \ref{copo1} we conclude that there exists a polynomial $w\in \C[z]$ of degree one such that either \be \label{krot} w\circ L_1= cz^r+\frac{1}{cz^r}\ee for some 
$c\in \C,$ or
\be \label{po} w\circ L_1=\bar L_1\circ 1/z.\ee

In the first case, it follows from the equalities 
\be \label{xer} L=(A\circ w^{-1})\circ (w\circ L_1)\ee and \eqref{krot} 
by Lemma \ref{uuii} that $\vert c\vert =1$ implying that   $w\circ L_1\in {\cal L_{\R}}[z]$. 
Now equality \eqref{xer} implies by Lemma \ref{l0} that $A\circ w^{-1}\in \R[z]$.
Thus,  we can set $v=w.$
 
Consider  the second case. Let $w=az+b,$ $a,b\in \C$, and $L_1=\sum_{-n}^nc_iz^i$, $c_i\in \C$. Then \eqref{po} implies the equalities 
$$  \bar c_{-i}=ac_i,\ \ \ 0 <\vert i \vert \leq n,$$ and therefore the equalities $$c_{-i}= \overline{a  c_i}=\bar a a c_{-i}.$$ Taking $c_{-i}\neq 0$, we conclude that 
$a\bar a=1$ or equivalently that $\vert a\vert =1.$ 
Setting now $v=\lambda z+\mu,$ where $\lambda$ satisfies $\lambda^2=a$ and $\mu=\overline{\lambda c_0}$, one can see easily that $v\circ L_1\in {\cal L_{\R}}[z]$. 
Indeed, the free term of $v\circ L_1$ is $\lambda c_{0}+\overline{\lambda  c_0}$ and therefore is real. For other terms, taking into account that  $\lambda\bar \lambda =1,$ 
we have:
$$\overline{\lambda c_{-i}}=\bar \lambda ac_{i}=\bar \lambda \lambda^2c_{i}=\lambda c_{i}, \ \ \ 0 <\vert i \vert \leq n.$$
Finally, Lemma \ref{l0} implies as above that $A\circ v^{-1}\in \R[z]$.
\qed

\bc \label{c1} Let $L=P\circ W$ be a decomposition of $L\in  {\cal L_{\R}}[z]\setminus \R$  into a composition
of $P,W\in \C(z)$. Then there exists a rational function $v\in \C(z)$ of degree one
such that $P\circ v^{-1}\in \R(x)$ and $v\circ W\in {\cal L_{\R}}(z)$.
\ec   
\pr  Arguing as in the proofs of Lemma \ref{l1} and Lemma \ref{l2} we see that 
there exists  a rational function $\mu\in \C(z)$ of degree one
such that either equality 
\eqref{iuy} holds or $P\circ \mu $ is a polynomial and $\mu^{-1}\circ W$ is a Laurent polynomial. In the first case, since 
$\mu^{-1}\circ W$ is contained in ${\cal L_{\R}}(z)$, it follows from Lemma \ref{l0} that  
$P\circ \mu \in \R(x),$ so we can set $v=\mu.$ In the second case the statement 
follows from Theorem \ref{2.1}

\section{Double decompositions in $\R_t[\theta]$ and in $\C[z,1/z]$}
For a rational function $P\in \C(z)$, 
two decompositions $P=A\circ B$ and $P=\widetilde A\circ \widetilde B$, where $A,B,\widetilde A,\widetilde B\in \C(z)$, 
are called equivalent if there exists a function $\mu\in \C(z)$ of degree one such that 
\be \la{pes} \widetilde A=A\circ \mu, \ \ \ \widetilde B=\mu^{-1}\circ B.\ee  Notice that if both $\widetilde A$ and $A$ (or $\widetilde B$ and $B$) are polynomials, then $\mu$ also is a polynomial. 
In particular, this is the case 
for most of the equivalences considered below. In case if we consider rational functions defined over an arbitrary field, the definition above is modified in an obvious way  (below we are only interested in the cases where the ground field is $\C$ or 
$\R$).
Abusing of notation we  will use for equivalent decompositions of rational functions the same symbol $\sim$ as for equivalent decompositions of trigonometric polynomials or polynomials.

We start from recalling some basic facts about polynomial solutions of the equation
\be \label{tor} A\circ C=B\circ D.\ee
The proposition below reduces a description of solutions of \eqref{tor} to the case where degrees of $A$ and $B$ as well as of $C$ and $D$ are coprime
(\cite{en}).

\bp \la{eng}
Suppose $A,B,C,D\in \C[z]\setminus \C$ satisfy \eqref{tor}. Then there exist 
$U, V, \tilde A, \tilde C, \tilde B, \tilde D \in \C[z], $  where
$$\deg U=\GCD(\deg A,\deg B),  \ \ \ \deg V=\GCD(\deg C,\deg D),$$
such that
$$A=U\circ \tilde A, \ \  B=U\circ \tilde B, \ \ C=\tilde C\circ V, \ \  D=\tilde D\circ V,$$
and $$\tilde A\circ \tilde C=\tilde B\circ \tilde D. \ \ \  \ \ \ \  \ \ \Box$$
\ep

\noindent In fact, under an appropriate restriction, Proposition \ref{eng} remains true if to assume that coefficients of polynomials  $A,B,C,D$ as well as of $U, V, \tilde A, \tilde C, \tilde B, \tilde D$ belong to an arbitrary field (see \cite{sch}, Chapter 1,  Theorem 5). In particular, Proposition \ref{eng} remains true if
the ground field is $\mathbb R.$

The following result obtained by Ritt \cite{r2} describes solutions of \eqref{tor} in the case where 
the equalities 
\be \label{qwer} \GCD(\deg A,\deg B)=1, \ \ \ \GCD(\deg C,\deg D)=1 \ee
hold, and is known under the name ``the second Ritt theorem''.

\bt \la{ritt} 
Suppose $A,B,C,D\in \C[z]\setminus \C$ satisfy \eqref{tor} and \eqref{qwer}.
Then 
there exist $U,\widetilde A, \widetilde B, \tilde C,\tilde D, W\in \C[z]$, where $\deg U=\deg W=1,$ such 
that $$ A=U \circ \widetilde A, \ \ \ \ B=U \circ \widetilde B, \ \ \ \ 
C=\tilde C \circ W, \ \ \ \ D=\tilde D \circ W,\ \ \ \ \widetilde A\circ \tilde C=\widetilde B\circ \tilde D$$
and, up to a possible replacement of $A$ by $B$ and $C$ by $D$, one of the following conditions holds:
$$\widetilde A\circ \tilde C\sim z^n \circ z^rR(z^n),  \ \ \ \ \ \ 
\widetilde B\circ \tilde D\sim  z^rR^n(z) \circ z^n,\leqno 1) $$ 
where $R\in \C[z]$, $r\geq 0,$ $n\geq 1,$  and 
$\GCD(n,r)=1;$
\vskip 0.01cm 
$$\widetilde A\circ \tilde C\sim T_n \circ T_m, \ \ \ \ \ \ \widetilde B\circ \tilde D\sim T_m \circ T_n,\leqno 2)$$ 
where $T_n,T_m$ are the  Chebyshev polynomials, $m,n\geq 1,$ and $\GCD(n,m)=1.$ \qed 

\et

Again, this theorem remains true if to assume that coefficients of all polynomials involved are real and, under an appropriate modification, even belong to an arbitrary field (see \cite{za} and \cite{sch}, Chapter 1, Theorem 8).

\vskip 0.2cm

Recall now the main result of the decomposition theory of Laurent polynomials (see \cite{pak})
concerning solutions of the equation 
\be \label{lau} P_1\circ W_1=P_2\circ W_2, \ee 
where $P_1,P_2\in \C[z]$ and $W_1,$ $W_2\in \C[z,1/z],$ using the  
notation of \cite{paen} (Theorem 3.1). 
Notice that the main result of \cite{paen} (Theorem A) also  may be used 
for a proof of Theorem 1.1. However, the approach based on the results of Section 2 is more general and 
may be used for a solution of other problems related to decompositions of trigonometric polynomials.

Set 
$$U_n=\frac{1}{2} \left(z^n+\frac{1}{z^n}\right), \ \ \  V_n=\frac{1}{2i} \left(z^n-\frac{1}{z^n}\right).$$
It is easy to see that the equalities 
$$\cos n\theta=T_n(\cos \theta), \ \ \ \sin n \theta=\frac{1}{n}T^{\prime}_n(\cos \theta)\sin \theta$$ and 
$$T_n\circ \frac{1}{2}\left(x+\frac{1}{x}\right)=\frac{1}{2}\left(x^n+\frac{1}{x^n}\right)$$
imply that 
$$\label{cos} U_n=\phi(\cos n\theta), \ \  \ V_n=\phi(\sin n\theta).$$  
Furthermore, if $c=\cos a+i\sin a,$ where $a\in \R$, then the equalities  
$$\cos (\theta+a)=\cos \theta \cos a-\sin \theta \sin a, \ \ \ \sin(\theta+a)=\sin \theta \cos a+\cos \theta \sin a,$$
imply that  
\be \label{cos1} U_n\circ (cz)=\phi(\cos (n(\theta+a))), \ \  \ V_n\circ(cz)=\phi(\sin (n(\theta+a))).\ee

\bt \la{irrr} 
Let $P_1,P_2\in \C[z]\setminus \C$ and $W_1,$ $W_2\in \C[z,1/z]\setminus \C$ 
satisfy \eqref{lau}.
Then 
there exist $F,$ $\widetilde P_1,$ $\widetilde P_2\in \C[z]$ and 
$W,$ $\tilde W_1,$ $\tilde W_2\in \C[z,1/z]$ such 
that
$$ P_1=F \circ \widetilde P_1, \ \ \ \ P_2=F \circ \widetilde P_2, \ \ \ \ 
W_1=\tilde W_1 \circ W, \ \ \ \ W_2=\tilde W_2 \circ W,\ \ \ \ \widetilde P_1\circ \tilde W_1=\widetilde P_2\circ \tilde W_2$$
and, up to a possible replacement of $P_1$ by $P_2$ and $W_1$ by $W_2$, one of the following conditions holds:
$$\widetilde P_1\circ \tilde W_1\sim z^n \circ z^rR(z^n),  \ \ \ \ \ \ 
\widetilde P_2\circ \tilde W_2\sim  z^rR^n(z) \circ z^n,\leqno 1) $$ 
where $R\in \C[z]$, $r\geq 0,$ $n\geq 1,$  and 
$\GCD(n,r)=1;$
\vskip 0.01cm 
$$\widetilde P_1\circ \tilde W_1\sim T_n \circ T_m, \ \ \ \ \ \ \widetilde P_2\circ \tilde W_2\sim T_m \circ T_n,\leqno 2)$$ 
where $T_n,T_m$ are the  Chebyshev polynomials, $m,n\geq 1,$ and $\GCD(n,m)=1;$ 
\vskip 0.01cm
$$\widetilde P_1\circ \tilde W_1\sim  z^2 \circ U_1
S(V_1),  \ \ \ \ \ \ \widetilde P_2\circ \tilde W_2\sim  (1-z^2)\,S^2 \circ V_1,\leqno 3)$$
where $S\in \C[z]$;
\vskip 0.01cm
$$\widetilde P_1\circ \tilde W_1\sim -T_{nl} \circ U_m(\v z), \ \ \ \ \ \ 
\widetilde P_2\circ \tilde W_2\sim T_{ml} \circ  U_n,\leqno 4)$$    
where $T_{nl},T_{ml}$ are the Chebyshev polynomials, $m,n\geq 1,$ $l>1$, $\varepsilon^{nlm}=-1$,
and $\GCD(n,m)=1;$ 
\vskip -0.2cm 
$$\widetilde P_1\circ \tilde W_1\sim (z^2-1)^3\circ \left(
\frac{i}{\sqrt{3}}\,V_2+\frac{2\sqrt{2}}{\sqrt{3}}\,U_1\right), \leqno 5) $$
$$\ \ \ \ \ \ \widetilde P_2\circ \tilde W_2\sim (3z^4-4z^3)\circ  \left(
\frac{i}{3\sqrt{2}}\,V_3+U_2
+\frac{i}{\sqrt{2}}\,V_1+\frac{2}{3} \right). \qed$$ 

\et

\vskip 0.2cm

\noindent Notice that if $W_1,W_2$ are polynomials, then $W$ also is a polynomial and  either 1) or 2) 
holds, in correspondence with Proposition \ref{eng} and Theorem \ref{ritt}.

\vskip 0.3cm

\noindent{\it Proof of Theorem 1.1.} Let $P_1, P_2 \in \R[x]$ and  $w_1,w_2\in \R_t[\theta]$ satisfy equation \eqref{maii}.  Assume first that there exist $w\in \R_t[\theta]$ and $\widehat W_1, \widehat W_2 \in \R[x]$ such that 
the equalities 
\be \label{vot} w_1=\widehat W_1\circ w, \ \ \ w_2=\widehat W_2\circ w\ee hold.
Then equality \eqref{maii} implies the equality 
$$P_1\circ \widehat W_1=P_1\circ \widehat W_1,$$ and it is easy to see using the real versions of Proposition \ref{eng} and Theorem \ref{ritt} that 
either the case 1, a) or the case 1, b)  of Theorem 1.1 has the place. 

Assume now that such $w$ and $\widehat W_1,$ $\widehat W_2$ do not exist. 
Set 
$$p=P_1\circ w_1= P_2\circ w_2, \ \  \ L=\phi(p), \ \ \  W_1=\phi(w_1), \ \ \ W_2=\phi(w_2),$$ and apply Theorem \ref{irrr} to equality \eqref{lau}.
Observe that our assumption implies that  neither the first nor the second case provided by Theorem \ref{irrr} may have the place. Indeed, since $L$ is a proper Laurent polynomial, if one of these cases holds, then the function $W$ also  is  a proper Laurent polynomial. Therefore,     
applying Theorem \ref{2.1} to the equality $W_1=\tilde W_1 \circ W$, we conclude that 
there exists a polynomial $v\in \C[z]$ of degree one
such that $\tilde W_1\circ v^{-1}\in \R[x]$ and $v\circ W\in {\cal L_{\R}}[z]$.
Furthermore, applying Lemma \ref{l0} to the equality $$W_2=(\tilde W_2 \circ v^{-1})\circ (v \circ W),$$ we conclude 
that $\tilde W_2\circ v^{-1}\in \R[x]$ implying that \eqref{vot} holds for  
$$\widehat W_1= \tilde W_1\circ v^{-1}, \ \ \ \widehat W_2= \tilde W_2\circ v^{-1}, \ \ \ w=\phi^{-1}(v\circ W).$$ 

Consider now one by one all the other cases possible  by Theorem \ref{irrr}.  
If holds 3), then 
there exist $\mu_1,\mu_2\in \C[z]$ of degree one and $S\in \C[z]$ such that 
\be  \label{koto} P_1=  F\circ z^2 \circ \mu_1, \ \ \   W_1 =\mu_1^{-1} \circ U_1S(V_1) \circ W,\ee 
and
\be   \label{kot} P_2=  F\circ (1-z^2)\,S^2\circ \mu_2 , \ \ \   W_2 =\mu_2^{-1}  \circ  V_1\circ W,\ee
for some $F\in \C[z]$ and $W\in {\cal L}[z].$
Furthermore, it follows from Lemma \ref{l1} that $W$ necessary has the form $W=cz^k,$ $c\in \C\setminus \{0\}.$ 

Let $\alpha$ be the leading coefficient of the polynomial $F$, and $d=\deg F.$ Setting $\mu_1=\alpha_1z+\beta_1,$ where $\alpha_1,\beta_1\in \C,$ we see that  the coefficients of
$z^{2d}$ and $z^{2d-1}$ of the polynomial $P_1$ are $c_{2d}=\alpha\alpha_1^{2d}$ and 
$c_{2d-1}=\alpha\alpha_1^{2d-1}\beta_1 2d$. Therefore, since 
$P_1\in \R[x]$, the number $$\frac{\beta_1}{\alpha_1}= \frac{c_{2d-1}}{2dc_{2d-1}}$$  is real and hence
$\mu_1=\alpha_1\tilde{\mu},$ where $\tilde{\mu}=z+(\beta_1/\alpha_1)\in \R[z].$ Thus, 
changing $\mu_1$ to $\tilde \mu$,  $F$ to $F\circ (\alpha_1^2 z)$, and $S$ to $(1/\alpha_1^2)$,
without loss of generality we may assume that  $\mu_1\in \R[x]$.  Since $\bar P_1=P_1$, this implies that $F\in \R[x]$.

Further, if $\mu_2^{-1}=\alpha_2 z+\beta_2 ,$ where $\alpha_2,\beta_2\in \C,$ then, since $W_2$ is contained in ${\cal L_{\R}}[z]$, the second equality in \eqref{kot} implies that $\beta_2\in \R$ and, by Lemma \ref{uuii}, that 
$\alpha_2\in \R$ and $\bar c=1/c$. Therefore, $\mu_2\in \R[x]$.
Furthermore, since  $\bar c=1/c$ and  $\mu_1\in \R[x]$, it follows from  $W_1\in{\cal L_{\R}}[z]$ that $S\in \R[x]$. 
 Finally, since $\vert c\vert =1,$ 
there exists $a\in \R$ such that $c=\cos a +i\sin a$, implying by \eqref{cos1} that 
$$w_1=\mu_1\circ \cos (k\theta +b)S(\sin(k\theta +b)),  \ \ \ w_2=\mu_2\circ \sin(k\theta +b),$$ 
where $b=ka.$ 
Thus, equalities  \eqref{koto} and \eqref{kot} lead to the case 2, a).

Consider now case 4).
In this case 
there exist $\mu_1,\mu_2\in \C[z]$ of degree one and $F\in \C[z]$ such that 
\be  \label{koto+} P_1=  F\circ -T_{nl} \circ \mu_1, \ \ \   W_1 =\mu_1^{-1} \circ U_m(\varepsilon z)  \circ W,\ee 
and
\be   \label{kot+} P_2=  F\circ T_{ml}\circ \mu_2, \ \ \   W_2 =\mu_2^{-1}  \circ  U_n\circ W,\ee
where  $\varepsilon^{nlm}=-1$ and $W=cz^k,$ $c\in \C\setminus \{0\}.$ As above, the second equality in \eqref{kot+} implies that 
$\bar c=1/c$ and $\mu_2\in \R[x].$ Then, using  $\mu_2\in \R[x]$ we see that the first equality in \eqref{kot+} implies that  $F\in \R[x]$, and using $\bar c=1/c$ we see that the second equality in \eqref{koto+} implies 
that $\mu_1\in \R[x].$ Therefore, taking into account formulas \eqref{cos1}, we conclude that  equalities \eqref{koto+} and \eqref{kot+} lead 
to the case 2, b).

Let us show finally that the case 5) 
cannot have a place. Assume the inverse. Then 
$$ W_1=\mu \circ  \left(
\frac{i}{\sqrt{3}}\,V_2+\frac{2\sqrt{2}}{\sqrt{3}}\,U_1\right)\circ (cz^k)=\mu \circ  \left(
\frac{1}{2\sqrt{3}}\left(z^2-\frac{1}{z^2}\right)+\frac{\sqrt{2}}{\sqrt{3}}\left(z+\frac{1}{z}\right)\right)\circ (cz^k),$$
where $\mu=\alpha z+\beta,$  $\alpha,\beta, c\in \C,$ and $\alpha\neq 0,$ $c\neq 0.$
Since  $W_1\in {\cal L_{\R}}[z]$, 
this
implies that
$$\bar \alpha \bar c^2=-\alpha/c^2, \ \ \ \bar \alpha \bar c=\alpha/c,$$ 
and dividing the first  equality by the second one we obtain the equality
$\bar c c=-1$ which is impossible.     \qed


\begin{thebibliography}{9}


\bibitem{abc} A. Alvarez, J. L. Bravo, C. Christopher, {\it On the trigonometric moment problem}, 
Ergodic Theory and Dynamical Systems,  34 (2014), no. 1, 1-20. 



\bibitem {bpy} M. Briskin, F. Pakovich, Y. Yomdin, \textit{Algebraic Geometry of the Center-Focus problem for Abel Differential Equation}, 
Ergodic Theory and Dynam. Sys., to appear. 



\bibitem {bry} M. Briskin, N. Roytvarf, Y. Yomdin,
{\it Center conditions at infinity for Abel differential equation}, Ann. of Math.,  172  (2010),  no. 1, 437-483.



\bibitem {cher} L. Cherkas,  {\it Number of limit cycles of an autonomous
second-order system}, Differentsial'nye uravneniya  {\bf 12} (1976),
No.5, 944-946




\bibitem {ga1} A. Cima, A. Gasull, F. Mañosas, 
{\it A simple solution of some composition conjectures for Abel equations,} J. Math. Anal. Appl. 398 (2013), 477-486.



\bibitem {en} H. Engstrom, \textit{Polynomial substitutions,} Amer. J. Math. 63, 249-255 (1941).





\bibitem{ggl} J. Gin\'e, M. Grau, J. Llibre, {\it Universal centres and composition conditions}, 
Proc. London Math. Soc.,  106 (2013), no. 3, 481-507. 






\bibitem{mp} M.\,Muzychuk, F.\,Pakovich, {\it 
Jordan-Holder theorem for imprimitivity systems and maximal decompositions of rational functions},  Proc. Lond. Math. Soc. (3)  102  (2011),  no. 1, 1-24.





\bibitem {pm} F. Pakovich, M. Muzychuk, {\it Solution of the polynomial moment problem}, Proc. Lond. Math. Soc.,  99  (2009),  no. 3, 633-657.



\bibitem{pak} F. Pakovich, \textit{Prime and composite Laurent polynomials},  Bull. Sci. Math.,  133  (2009),  no. 7, 693-732.



\bibitem{pakk} F. Pakovich, \textit{Generalized ``second Ritt theorem" and explicit solution of the polynomial moment problem}, Compositio Math., 149 (2013), 705-728.


\bibitem{ppre} F. Pakovich, \textit{On rational functions orthogonal to all powers of a given rational function on a curve}, 
 Mosc. Math. J. 13 (2013), no. 4, 693-731. 




\bibitem{ppz} F. Pakovich, C. Pech, A. Zvonkin, \textit{ Laurent polynomial moment problem: a case study}, 
Contemp. Math., 532 (2010), 177-194.


\bibitem{paen} F. Pakovich, \textit{On the equation P(f)=Q(g), where P,Q are polynomials and f,g are entire functions}, Amer. Journal of Math., 132 (2010), no. 6, 1591-1607.


\bibitem{pakov} F. Pakovich, \textit{ Weak and strong composition conditions for the Abel differential equation},   Bull. Sci. Math. 138 (2014), no. 8, 993-998.





\bibitem{r2} J. F. Ritt. \textit{Prime and composite polynomials,} Trans. Amer. Math. Soc., 23, 1922, 51-66.



\bibitem{r} J. F. Ritt. \textit{Permutable rational functions,} Trans. Amer. Math. Soc. 25 (1923), 399-448. 



\bibitem {sch} A. Schinzel, \textit{Polynomials with special regard to
reducibility}, Encyclopedia of Mathematics and Its Applications
77, Cambridge University Press, 2000.

\bibitem {za} U. Zannier, 
\textit{Ritt's second theorem in arbitrary characteristic,}
J. Reine Angew. Math. 445, 175-203 (1993)


 
\end{thebibliography}
\end{document}